
\hbadness = 1000	

\vbadness = 10000	

\hfuzz=10pt	

\input amsppt.sty
\output={\plainoutput}
\magnification \magstep1
\parskip=.1in
\hsize=6.5 true in

\baselineskip=18pt  

\define\theorem#1.#2{{\smc #1.#2 Theorem.}$\ $}
\define\corollary#1.#2{{\smc #1.#2 Corollary.}$\ $}
\define\lemma#1.#2{{\smc #1.#2 Lemma.}$\ $}
\define\proof{{\smc Proof.}$\ $}

\define\endproof{$ \qquad \square $}

\def \ss#1{\lower3pt\hbox{${\scriptstyle #1}$}}

\define\Hi{H^\infty}

\define\s{\Cal}

\centerline{{\smc Analytic Disks in Fibers over the Unit Ball of
    a Banach Space}}
\vskip .2in
\centerline{{\smc B.J. Cole,*\footnote" "{*Supported by NSF grant
    \#DMS86-12012}
T.W. Gamelin**\footnote" "{**Supported by NSF grant \#DMS88-01776}
and W.B. Johnson***\footnote" "{***Supported by NSF grant
\#DMS90-03550}}}

\noindent
Mathematics Department
\newline
Brown University
\newline
Providence, RI 02912

\noindent
Mathematics Department
\newline
UCLA
\newline
Los Angeles, CA 90024

\noindent
Mathematics Department
\newline
Texas A \& M University
\newline
College Station, TX 77843

\vskip .3in
{\bf Abstract.} \
We study biorthogonal sequences with special properties, such as
weak or weak-star convergence to $0$, and obtain an extension of the
Josefson-Nissenzweig theorem.  This result is applied to embed analytic
disks in the fiber over $0$ of the spectrum of $\Hi (B)$, the algebra of
bounded analytic functions on the unit ball $B$ of an arbitrary infinite
dimensional Banach space.  Various other embedding theorems are
obtained. For instance, if the Banach space is superreflexive, then
the unit ball of a Hilbert space of uncountable dimension can be
embedded analytically in the fiber over $0$ via an embedding which is
uniformly bicontinuous with respect to the Gleason metric.


\newpage

{\bf 1. Introduction}. \ Fix an infinite dimensional complex Banach
space $\s X$, with open unit ball $B$.  We are interested in studying
the uniform algebra $\Hi (B)$ of bounded analytic functions on $B$,
and its spectrum ${\s M} = {\s M} (B)$ consisting of the nonzero
complex-valued homomorphisms of $\Hi (B)$.  The spectrum $\s M$ is
fibered in a natural way over the closed unit ball $\bar B^{**}$ of
the bidual ${\s X}^{**}$ of $\s X$.  The projection of $\s M$ onto
$\bar B^{**}$ is obtained by simply restricting $\varphi \in \s M$ to
$\s X^*$, regarded as a subspace of $\Hi (B)$.  As a straightforward
application of the Josefson-Nissenzweig theorem [Jo,Ni], it is shown
in [ACG] that each fiber consists of more than one point and is in
fact quite large.  Our aim here is to prove a sharpened form of the
Josefson-Nissenzweig theorem, and to use this to embed analytic disks
in the fiber over $0$.  This stands in contrast to the situation in
finite dimensional Banach spaces, where one expects (and can prove under
certain hypotheses) that the natural projection is one-to-one over the
open unit ball $B = B^{**}$ and implements a homeomorphism of $B$ and an
open subset of $\s M$.

The Josefson-Nissenzweig theorem asserts that in any infinite
dimensional dual Banach space $\s Z$, there is a sequence $\{ z_j \}$
converging weak-star to $0$, such that $|| z_j || = 1$.  The accessory
condition we require is that the distance from $z_j$ to the linear
span of the preceding $z_i$'s tends to $1$ as $j \rightarrow \infty$.
Actually we prove that each $z_j$ can be chosen to have unit distance
from the linear span of the remaining $z_i$'s, and in fact the $z_j$'s
can be taken as part of a unit biorthogonal system, defined in Section
2.

In Sections 2 and 3 we establish the existence of unit
biorthogonal systems having accessory properties involving
weak or weak-star convergence to $0$.
In Section 4 we make some observations regarding infinite products in
a uniform algebra.  The embedding of analytic disks is accomplished in
Section 5. In Section 6 we give conditions under which the unit ball
of certain infinite dimensional Banach spaces $\s Y$ can be injected
analytically into the fiber over $0$.  For instance, when $\s X$ is
$c_0$, $\s Y$ can be $\ell ^ \infty$ and the injection isometric
with respect to the Gleason metric, while
if $\s X$ is superreflexive and infinite dimensional, $\s Y$ can be
a non-separable Hilbert space and the injection
uniformly bicontinuous. In Section 7 we indicate how certain of the
embedding results extend to other algebras of analytic functions
associated with $\s X$.

For background on analytic functions on Banach spaces see [Mu], and for
uniform algebras see [Ga].  When convenient we will regard the functions
in the algebra at hand to be continuous functions on the spectrum of
the algebra, via the Gelfand transform.  We denote the open unit disk
in the complex plane {\bf C} by $\Delta$.  The linear span of a set of
vectors $S$ will be denoted by $sp(S)$, and its closure by $ \overline
{sp} (S)$.  The distance from a vector $x$ to a subset $Y$ of $\s X$,
measured in the norm of $\s X$, will be denoted by $dist(x,Y)$.

\vskip .3in

{\bf 2. Unit biorthogonal systems.} \
Let  ${\s X}$  be an infinite dimensional Banach space,
with dual ${\s X}^*$.
A {\it unit biorthogonal system} for  ${\s X}$
consists of sequences  $\{x_j\}$   in  ${\s X}$  and
$\{x_k^*\}$  in  ${\s X}^*$  satisfying
$$
||x_j|| = 1 = ||x_k^*||,
\tag{2.1}
$$
$$
x_k^* (x_j) = \cases 1,\qquad  j = k, \\
0,\qquad j \neq k. \endcases
\tag{2.2}
$$
We will be interested in finding such systems with the additional property
that  $x_j$  converges weakly to  $0$  in ${\s X}$.  The standard
basis and dual basis for  $\ell^p, \ 1 < p < \infty$,  form such a system.
However,  $\ell^1$  cannot have such a system on account of the Schur
property: weakly convergent sequences in  $\ell^1$  converge in norm.
If a subspace  ${\s Y}$  of  ${\s X}$  has such a
system, then by extending the  $x_j^*$'s via the Hahn-Banach theorem we
obtain a unit biorthogonal system for ${\s X}$.

To construct unit biorthogonal systems we use the following lemma which
is due to Krasnoselskii, Krein and Milman [KKM].

\lemma 2.1
{\it Let ${\s X}_0$ and ${\s X}_1$ be finite-dimensional subspaces
of  ${\s X}$ such that the dimension of ${\s X}_1$ exceeds that
of ${\s X}_0$.
Then there is a vector  $x \in {\s X}_1$
such that $||x|| = 1 = \text{dist} (x,{\s X}_0).$}

\proof
The proof is a simple application of Borsuk's theorem, that there
is no continuous antipodal map from one sphere to another of lower
dimension.
By perturbing the norm and taking limits, we can assume
that every $x \in {\s X}_1$ has a unique nearest point  $\phi(x) \in
{\s X}_0$.  We must find  $x \in {\s X}_1$ such that  $||x|| =
1$  and  $\phi (x) =0.$  If there were no such point, then  $x \rightarrow
\phi (x) / ||\phi (x)||$  would be an antipodal map  $(\phi(-x) = -
\phi(x))$  from the unit sphere of ${\s X}_1$  to the
unit sphere of ${\s X}_0$, contradicting Borsuk's theorem.
\endproof

Our basic method consists of repeated applications of the
following.

\lemma 2.2
{\it
Let $\{x_i\}_1^{n-1}$,$\{x_i^*\}_1^{n-1}$
be a unit biorthogonal system for ${\s X}$ of length $n-1$,
and let ${\s Z}$ be a
subspace of ${\s X}$ of dimension at least $2n-1$.  Then there exist
$x_n \in {\s Z}$ and $x_n^* \in {\s X}^*$ so that
$\{x_i\}_1^n$and $\{x_i^*\}_1^n$ form
a unit biorthogonal system of length $n$.}

\proof
The dimension of the space ${\s Z}_1$ of $z \in {\s Z}$ satisfying
$x_j^*(z) = 0$ for $1 \leq j \leq n-1$ is at least $n$.
By Lemma 2.1, there is  $x_n \in {\s Z}_1$   such that
$$
||x_n || = 1 = \text{dist} (x_n, \ \text{sp}\{x_1,\ldots,x_{n-1}\}).
$$
By the Hahn-Banach theorem, there is  $x_n^* \in {\s X}^*$  such that
$||x_n^*|| = 1 = x_n^* (x_n),$  and  $x_n^* (x_j) = 0$  for  $1  \leq j
< n-1$.
\endproof

\lemma 2.3
{\it
Let $\{y_i\}$ be a basic sequence in ${\s X}$. Then there are a
block basis $\{x_j\}$  of  $\{y_i\}$  and functionals $\{x_k^*\}$
in ${\s X}^*$ which form a unit biorthogonal system for  ${\s X}$.}

\proof
Suppose that  $x_1,\ldots,x_{n-1} \in {\s X}$, $x_1^*,
\ldots,x_{n-1}^* \in {\s X}^*$  and integers  $1 \leq m_1 < m_2 <
\cdots < m_{n-1}$  have been chosen to satisfy (2.1) and (2.2), so that
each  $x_j$  is a linear combination of the  $y_i$'s  for
$m_{j-1} < i \leq m_j$.  Choose $m_n = m_{n-1} + 2n - 1,$
and apply Lemma 2.2 to the linear span
of the $y_i$'s for $m_{n-1} < i \leq m_n$.
The sequences constructed inductively in this way form a unit
biorthogonal system.
\endproof

\theorem 2.4
{\it  For any infinite dimensional Banach space ${\s X}$,
there is a unit biorthogonal system
$\{x_j\}$, $\{x_k^*\}$  such that  $0$  is in the weak
closure of the  $x_j$'s.}

\proof
The  $x_j$'s and  $x_k^*$'s  are constructed by an induction process
similar to the procedure in the proof of Lemma 2.3, but with the
following twist.
By Dvoretzky's spherical sections theorem
(cf. [MS]), there is for any $\varepsilon > 0$
a subspace of ${\s X}$ of
arbitrarily large dimension $N$ which is $(1+\varepsilon )$-isomorphic
to  $\ell_N^2$.
Now we repeatedly invoke Lemma 2.2 while arranging that successively
larger blocks of the $x_j$'s lie in almost-euclidean
subspaces ${\s M}_N$ of ${\s X}$. The $x_j$'s need not be orthogonal
with respect to the euclidean structure of ${\s M}_N,$ but they are
`almost' orthonormal, and the `almost' can be expressed
quantitatively in terms of $\varepsilon$, so that for
$\varepsilon = \varepsilon _N$ sufficiently small
we obtain an estimate of the form
$$
\sum _{x_j \in {\s M}_N} |x^* (x_j)|^2 \leq C ||x^*||^2, \qquad
x^* \in {\s X}^*,
$$
where $C$ is independent of $N$.
Since the number of summands grows larger with $N$,
an increasingly large proportion of the values  $x^*(x_j)$
are small.  It follows that  $0$  is a weak cluster point of
the sequence  $\{x_j\}$.
\endproof

\theorem 2.5
{\it  Let  ${\s X}$  be an infinite dimensional Banach space
such that $\ell^1$  does not embed into
${\s X}$.  Then there are sequences  $\{x_j\}$  in  ${\s X}$
and  $\{x_k^*\}$ in ${\s X}^*$
forming a unit biorthogonal system, such that $x_j$
converges weakly to  $0$.}

\proof
Choose sequences $\{x_j\}$, $\{x_k^*\}$ as in Theorem 2.4.
Replacing $\s X $ by the closed linear span of the $x_j$'s, we can
assume that $\s X$ is separable, so that the weak-star topology on
the closed unit ball ${\bar B}^*$ of ${\s X}^*$ is a compact metric
topology.  By a theorem of Bourgain, Fremlin and Talagrand ([BFT];
see also [Di2, p.216]) the functions on ${\bar B}^*$
in the first Baire class, endowed with the topology of pointwise
convergence, form an angelic space $\Omega$.  Since $\s X$ contains no
copies of ${\ell}^1$, Rosenthal's theorem (cf. [Di2,
Chapter XI]) shows that the sequence $\{ x_j \}$, viewed as a
sequence of functions on ${\bar B}^*$, is sequentially
precompact in $\Omega$.  Since $0$ lies in the
pointwise closure of the sequence (by Theorem 2.4),
the angelic property provides
us with a subsequence which converges weakly to $0$.  Passing to
subsequences then, we obtain the desired unit biorthogonal system.
\endproof

Odell, Rosenthal, and Schlumprecht [ORS] have obtained results which
show in particular that in every infinite dimensional Banach space which
does not contain a copy of $\ell^1$ there is a basic sequence $\{y_i\}$
such that if $x_j$ is in the unit ball of the linear span of the $y_i$'s
for $2^j < i \leq 2^{j+1}$, then $\{x_j\}$ converges weakly to $0$.
This result can be used instead of the theorem from [BFT] to prove
Theorem 2.5.

The following result will be used for the embedding
of analytic disks in Section 5.

\theorem 2.6
{\it  If ${\s X}$ is any infinite dimensional Banach space,
there is a unit biorthogonal system $\{ x_j^{**}\}$,
$\{ x_k^{***}\}$ for the bidual
${\s X}^{**}$ of $\s X$ such that $\{ x_j^{**}\}$
converges weakly to $0$ in ${\s X}^{**}$.}

\proof
In view of Theorem 2.5, it suffices to find an infinite
dimensional subspace ${\s Y}$ of ${\s X} ^{**}$
into which $\ell ^1$ does not embed.
If $\ell ^1$ does not embed in ${\s X}$,
then take ${\s Y} = {\s X}$.
Suppose on the other hand that $\ell ^1$ embeds in $\s X$.
Then $(\ell ^ \infty )^*$ embeds in ${\s X}^{**}$.
Now $L^1[0,1]$ embeds in $(\ell ^ \infty )^*$, for instance in
$L^1(\nu)$ for an appropriate measure $\nu$ on the Stone-\v Cech
compactification of the integers, and
$\ell ^2$ embeds in $L^1[0,1]$, for instance as lacunary Fourier
series.
In this case, take ${\s Y}$ to be the embedded image of
$\ell ^2$ in ${\s X}^{**}$.
\endproof

\vskip .3in

{\bf 3. Unit biorthogonal systems in a dual Banach space.} \ In this
section we aim to establish the existence of biorthogonal systems for a
dual Banach space which converge weak-star to $0$.  The proof will be
broken into cases, according to the following lemma.

\lemma 3.1
{\it Let  ${\s Z}$  be a dual Banach space.  Then
at least one of the following conditions holds:
\roster
\item "{(i)}" There is an embedding of $\ell^1$ into ${\s Z}$ such
that the image of the standard basis converges weak-star to $0$.
\item "{(ii)}" There is an infinite dimensional subspace ${\s Y}$
of ${\s Z}$ into which $\ell^1$ does not embed.
\endroster
}

\proof
Suppose that $\ell^1$ embeds in ${\s Z}$
and (i) fails, that is,
there is no embedding for which the image of the standard basis
tends weak-star to $0$.  By [HJ] (or see [Di2, p.219]), $\ell^1$ embeds
in the predual $\s X$ of $\s Z$.  From a theorem of Pe\l czy\'nski
([Pe, Theorem 3.4]; [Di2, p.213]) it follows that $L^1[0,1]$ embeds in
${\s X}^* = {\s Z}$.  Since $\ell^2$ embeds in $L^1[0,1]$, we can 
take ${\s Y}$ to be the embedded image of $\ell^2$ in ${\s Z}$.
\endproof

\theorem 3.2
{\it  Let  ${\s Z}$  be an infinite dimensional
dual Banach space. Then there is a unit biorthogonal system  $\{z_j\},$
$\{z_k^*\}$  for  ${\s Z}$  such that  $z_j$  converges weak-star to
$0$. Either the $z_j$'s can be chosen to converge weakly to $0$; or
else there is an embedding of  $\ell^1$  into  ${\s Z}$  such that
the image  $\{y_i\}$  of the standard basis for  $\ell^1$  converges
weak-star to  $0$, in which case the $z_j$'s can be chosen to be
a block basis of the $y_i$'s.}

\proof
Consider the two cases of Lemma 3.1. In case (ii)
the unit biorthogonal system is provided
by  Theorem 2.5.

Suppose that we are in case (i) of Lemma 3.1.
Applying Lemma 2.3 to the image  $\{y_i\}$  of the standard basis for
$\ell^1$, we obtain $z_j$'s and  $z_k^*$'s which form a unit
biorthogonal system, such that the $z_j$'s are a block basis
of the $y_q$'s. Then
each  $z_j$  has the form  $\sum a_{jq} y_q,$  and since the  $z_j$'s
have unit norm, we obtain a uniform bound
$$
\sum_q |a_{jq}| \leq c, \quad 1 \leq j < \infty.
$$
Hence if  $x \in {\s X}$,  the predual of  ${\s Z}$,  we obtain
$$
|z_j(x)| = |\sum a_{jq} y_q (x)| \leq c \underset{n_{j-1} < q \leq n_j}\to
\max |y_q (x)|.
$$
Since  $y_q(x)$  tends to  $0$,  so does  $z_j(x),$  and  $z_j$  converges
weak-star to  $0$.
\endproof

While we are focusing on complex Banach spaces, the results of
Sections 2 and 3 are all valid also for real Banach spaces.

\vskip .3in

{\bf 4. Infinite products in a uniform algebra.}$\ $In this section we
record some observations on the convergence of infinite products in a
uniform algebra.  For Blaschke products in this context, see [Gl].

Let $A$ be a
uniform algebra,  with spectrum  $M_A$.  We regard  $A$  as a
uniform algebra on  $M_A$. For $\psi \in M_A$ and $0<r<1$, let
$S_r(\psi)$ denote the hyperbolic ball of radius $r$ centered at $\psi$,
defined as the collection of $\varphi \in M_A$ such that
$|f(\varphi)| \leq r$ for all $f \in A$ satisfying $||f|| \leq 1$
and $f(\psi) = 0.$ Thus the Gleason part of $\psi$ is the union of
the $S_r$'s for $0<r<1.$  Harnack's estimate
for uniform algebras [Ga, Section VI.1] shows that
if $|f_j| \leq 1$ and
$f_j(\psi) \rightarrow 1,$ then $f_j$ tends uniformly to $1$ on
each hyperbolic ball $S_r(\psi)$. Applying this remark to partial
products, we are led to the following.

\lemma 4.1
{\it
Fix $g_1, g_2,\ldots \in A$ satisfying $|g_j| \leq 1$,
and consider the infinite
product $G$ defined wherever it converges by
$$
G(\varphi) = \prod_{j=1}^\infty g_j(\varphi).
$$
If the product converges at $\psi \in M_A$ to a nonzero
value, then it converges uniformly
on each hyperbolic ball $S_r(\psi)$, $0<r<1$.
}

We are interested in infinite products formed from elements
of $\Hi (B)$.  Before discussing convergence, we provide some
background. According to [DG, Theorem 5],
the functions in $\Hi (B)$ have natural extensions to the open
unit ball $B^{**}$ of ${\s X}^{**}$, so that $\Hi (B)$ is embedded
isometrically as a closed subalgebra of $\Hi (B^{**})$.  We can
thus regard $B^{**}$ as a subset of $\s M$, by identifying
$z \in B^{**}$ with the homomorphism which evaluates the natural
extension of $f \in \Hi (B)$ at $z$. From the Schwarz lemma it is
easy to see that the hyperbolic balls centered at $0$ meet $B^{**}$
in closed norm balls:
$$
S_r(0) \cap B^{**} = r \bar B^{**}, \qquad 0<r<1.
$$

Following [ACG], we define a function $g \in \Hi (B^{**})$ to be
{\it canonical } if $g$ is the natural extension to $B^{**}$ of
its restriction to $B$.  If $g \in \Hi (B^{**})$ is a pointwise
limit on $B^{**}$ of a bounded sequence of canonical functions
which converges uniformly on each $rB^{**}$,
then $g$ is canonical [ACG, Lemma 10.3].  Combining these
observations with Lemma 4.1,
we obtain the following.

\lemma 4.2
{\it Suppose $\{g_j\}$ is a sequence of functions in $\Hi (B)$
satisfying $|g_j| \leq 1$, and suppose the infinite product
$\prod g_j $ converges at $0$ to a nonzero value.
Then the infinite product converges uniformly on
$rB$ for each $r<1$ to $G \in \Hi (B)$.  Furthermore, the
canonical extensions to $B^{**}$ of the partial products converge
uniformly on each $rB^{**}$ to the canonical extension of $G$.}

Note that any subproduct of $G$ divides $G$ in $\Hi (B)$.
Consequently if $\varphi \in \s M $ satisfies
$\varphi (g_m) = 0 $ for some $m$, then
$\varphi (G) = 0$, because $g_m$  divides $G$.

\vskip .3in

{\bf 5. Embedding analytic disks in ${\s M}_0$.}$\ $
A sequence $S$ in $B$ is an {\it interpolating sequence for
$\Hi (B)$ } if the restriction of
$\Hi (B)$ to $S$ coincides with ${\ell}^{\infty}$.
We extend this definition to
sequences $S$ in $B^{**}$ by defining such a sequence to
be {\it interpolating for
$\Hi (B)$ } if the restriction to $S$ of the canonical functions in
$\Hi (B^{**})$ coincides with ${\ell}^{\infty}$.

The Blaschke products of interpolating sequences in the
open unit disk $\Delta$ in the complex plane can be used
to produce analytic disks in the fibers of $\s M (\Delta )$
over points of the boundary of $\Delta$.  For details, see [Ho].
We will use Blaschke products of interpolating sequences in $B^{**}$
to embed analytic disks in fibers over $B^{**}$.

\theorem 5.1
{\it Let  ${\s X}$  be an
infinite-dimensional Banach space.  Suppose
$\{z_k\}$ is a sequence in $B^{**}$ which
converges weak-star to $0$, such that
the distance from $z_k$ to the linear span of
$z_1,\ldots,z_{k-1}$ tends to 1 as $k \rightarrow \infty$.
Then passing to a subsequence we can find a sequence of
analytic disks $\lambda \rightarrow z_k(\lambda), \, \lambda
\in \Delta , \,
k\geq 1,$ in $B^{**}$ with $z_k(0) = z_k$, such that for each
$\lambda \in \Delta $, $\{z_k(\lambda )\}$ is an interpolating
sequence for $\Hi (B)$.  Furthermore, the correspondence
$(k,\lambda) \rightarrow z_k(\lambda)$ extends to an embedding
$$
\Psi : \beta (N) \times \Delta \rightarrow {\s M}
$$
such that
$$
\Psi ((\beta (N) \backslash N) \times \Delta) \subset {\s M}_0,
$$
and  $f\,{}_{{}^\circ} \Psi$  is analytic on each slice  $\{p\} \times
\Delta,$  for all  $f \in H^\infty (B)$  and  $p \in \beta (N)$.}

\proof
We can pass to a subsequence of the  $z_k$'s,  whenever
appropriate.  Thus we can assume that  $||z_k||$  converges very rapidly
to 1, and that in fact
$$
\text{dist} (z_k, \text{sp}\{z_1,\ldots,z_{k-1}\}) \rightarrow 1
$$
very rapidly.

Fix  $\delta > 0$  small.  For each  $z_k$,  choose  $r_k > 0$  so that
$$
r_k < ||z_k||,
\tag{5.1}
$$
$$
\frac{r_k(1 - ||z_k||)}{||z_k|| - r_k^2} > \delta .
\tag{5.2}
$$
The condition (5.2) is satisfied for  $r_k  = ||z_k||,$  so both (5.1)
and (5.2) are satisfied for  $r_k$  slightly less than  $||z_k||$.  Fix
$\varepsilon > 0$  small.  Choose $0 < s_k < 1$ such that
$\sum (1 - s_k) < \infty$, and such that
(passing to a subsequence if necessary)
$$
\bigg|\frac{r_j - \zeta }{1 - r_j \zeta }\bigg| > s_j ,
\qquad |\zeta | < \varepsilon .
\tag{5.3}
$$
The condition (5.3) means simply that the pseudohyperbolic disk centered at
$r_j$  with radius  $s_j$  is disjoint from the disk  $\{|\zeta| <
\varepsilon\}.$

Fix  $\beta > 0$  small.  We claim that, passing to a subsequence if
necessary, we can find  $L_k \in {\s X}^*$,  $k \geq 1,$  such that
$$
||L_k|| < 1,
\tag{5.4}
$$
$$
L_k (z_k) = r_k,
\tag{5.5}
$$
$$
L_k (z_j)  = 0, \qquad 1 \leq j < k,
\tag{5.6}
$$
$$
|L_k (z_j)| < \beta/ 2^k, \qquad j > k .
\tag{5.7}
$$
Indeed, suppose we have chosen, after discarding some  $z_j$'s  and
relabeling, functionals   $L_1,\ldots,L_{m-1}$  satisfying (5.4)  to
(5.7) for  $1 \leq j, \ k \leq m - 1.$  Since $z_j$ tends weak-star
to $0$, we can
also arrange after more discarding and relabeling that
$$
|L_j (z_m)| < \beta / 2^m, \qquad 1 \leq j \leq m - 1.
$$
Using the distance hypothesis, we can then apply the Hahn-Banach
theorem to find
$\Lambda \in
(sp\{z_i\}_{i=1}^m)^*$
such that  $||\Lambda|| < 1, \  \Lambda(z_m) > r_m,$
and  $\Lambda (z_k) = 0$  for
$1 \leq k < m$.
Since the restriction mapping from ${\s X}^*$ to
$(sp\{z_i\}_{i=1}^m)^*$ is a quotient mapping,
we obtain  $L_m \in {\s X}^*$  with the asserted properties.

Now define
$$
w_k(\lambda) = \frac{r_k - \lambda}{1 - \lambda r_k}
\frac{z_k}{r_k}, \qquad
|\lambda| \leq \delta .
$$
Then  $w_k (0) = z_k$.  Since the maximum of  $|r_k -
\lambda | / |1 - \lambda r_k|$  over the disk  $|\lambda| \leq \delta$
is attained at  $\lambda = - \delta$,  we obtain the estimate
$$
||w_k(\lambda)|| \leq \frac{r_k + \delta}{1 + \delta r_k}
\frac{||z_k||}{r_k},
\qquad |\lambda| \leq \delta .
$$
The condition (5.2) on  $r_k$  is equivalent to the right-hand side
of this estimate being less than $1$, so that we obtain
$$
||w_k(\lambda) || < 1, \qquad |\lambda| < \delta .
$$

We define a Blaschke product $G$ by
$$
G(z) = \prod_{j=1}^\infty \frac{r_j - L_j (z)}{1 - r_j L_j(z)} ,
\qquad z \in B^{**} .
$$
By Lemma 4.2 the product converges uniformly on  $rB^{**}$  for each
$0 < r < 1,$  and  $G \in H^\infty (B^{**})$  satisfies  $|G| < 1.$
Furthermore the function
$G$  is
canonical;
that is,  $G$  represents the canonical
extension to  $B^{**}$  of its restriction  $G|_B$  in  $H^\infty (B)$.

The factors of the Blaschke product  $G$  are chosen so that
$$
\frac{r_k - L_k (w_k (\lambda))}{1 - r_k L_k (w_k(\lambda))} = \lambda .
$$
Hence we can express
$$
G(w_k(\lambda)) = \lambda g_k(\lambda),
$$
where $g_k(\lambda)$ is the product over $j \neq k$ of
$$
\frac{r_j - L_j (w_k (\lambda))}{1 - r_j L_j (w_k(\lambda))} .
$$
>From (5.6) we see that this reduces to $r_j$ for $j > k$.
For $1 \leq j < k$ one computes that the factor has the form
$r_j + O(\Lambda _j(z_k))$, where the error estimate is uniform for
$|\lambda| < \delta$ and $j \neq k$.
The estimate
$$
\sum_{j=1}^{k-1} |L_j (z_k) | \leq \beta \bigg(1 + \frac{1}{2} + \cdots
+ \frac{1}{2^{k-1}}\bigg) < 2 \beta
$$
shows that by choosing  $\beta$  sufficiently small and the  $r_j$'s
sufficiently close to 1, we can arrange that $g_k(\lambda)$ is
close to $1$ for all
$|\lambda| \leq \delta$, uniformly in $k$. Hence if
$|\lambda| < \delta/2,$  there is a unique  $\zeta_k(\lambda)$  satisfying
$|\zeta_k (\lambda) | < \delta$  and
$$
G(w_k (\zeta_k(\lambda))) =\lambda, \qquad |\lambda| < \delta/2.
$$
Now we define $z_k(\lambda)$ to be the reparametrization of
the $w_k$'s given by
$$
z_k(\lambda) = w_k (\zeta_k (\lambda)), \qquad|\lambda| < \delta/2.
$$
Then $z_k(\lambda)$ is a multiple of $z_k$ which depends analytically
on $\lambda$ and satisfies
$$
\align
& ||z_k (\lambda) || < 1, \qquad |\lambda| < \delta/2, \\
& z_k (0) = z_k , \\
& G(z_k(\lambda)) = \lambda , \qquad |\lambda| < \delta /2.
\endalign
$$

We claim that for each fixed  $\lambda, \ |\lambda| < \delta/2$,  the
sequence  $\{z_k (\lambda)\}$  is an interpolating sequence for  $H^\infty
(B)$.  For this, it suffices to find for each subset  $J \subset N$  a
canonical function  $f_J \in H^\infty (B)$  such that  $|f_J| \leq 1$  on
$B$,  $|f_J(z_k(\lambda))| < 1/3$   for  $k \in J$,  and  $|f_J(z_k
(\lambda)) - 1 | < 1/3$  for  $k \in N \backslash J$.  A function  $f_J$
having these properties is the Blaschke subproduct
$$
f_J(z) = \prod_{k \in J} \frac{r_k - L_k (z)}{1 - r_k L_k(z)} .
$$
From
$$
f_J (w_k(\lambda)) = \lambda \prod_{j\in J, \, j \neq k} \
\frac{r_j - L_j (w_k (\lambda))}{1 - r_j L_j (w_k (\lambda))}
$$
we obtain
$$
|f_J (z_k(\lambda)) | \leq \delta, \qquad k \in J, \ |\lambda| \leq
\delta/2.
$$
On the other hand, the earlier estimates for $g_k(\lambda)$  are
easily modified to show that   $f_J (z_k(\lambda))$  is near 1 when  $k
\notin J.$  The claim is established.

Now let  $D$  be the open disk  $\{|\lambda| < \delta/2\}$  in the complex
plane.  It suffices to establish the theorem, with  $\Delta$  replaced by
$D$.  We define  $\Psi$  on  $N \times D$  by
$$
\Psi (k,\lambda) = z_k (\lambda), \qquad |\lambda| < \delta/2, \
1 \leq k < \infty,
$$
where we regard  $z_k(\lambda)$  as a point in  ${\s M}$,
i.e., we
identify  $z_k(\lambda)$  with the evaluation homomorphism.
Since  ${\s M}$  is compact,  $\Psi$   extends to a map
$$
\Psi : \beta (N) \times D \rightarrow {\s M},
$$
which is continuous on  $\beta(N)$  for each fixed  $\lambda \in D.$  For
all  $f \in H^\infty(B),$  $ f\, {}_{{}^\circ} \Psi$  is bounded and
analytic on each slice  $\{p\} \times D$  of  $\beta (N) \times D,$
hence equicontinuous on the slices  $\{p\} \times \{|\lambda| < r
\delta/2\}, \ r < 1.$  It follows that  $\Psi$  is jointly continuous.
Since  $\Psi(N \times \{\lambda\})$  is an interpolating sequence, for
each fixed  $\lambda$, $\Psi$  is one-to-one on  $\beta(N) \times
\{\lambda\}.$  On account of the identity
$$
G(\Psi(p,\lambda)) = \lambda , \qquad (p,\lambda) \in \beta (N) \times D,
$$
$\Psi$  is one-to-one on $\beta(N) \times D$. By shrinking $\delta$
slightly we obtain then that $\Psi$ is an embedding. Finally note that
$z_k(\lambda)$  converges weak-star to $0$  in  ${\s X}^{**},$  so
that  $\Psi$  maps  $(\beta(N) \backslash N) \times D$  into  ${\s
M}_0$.
\endproof

By Theorem 2.6, or Theorem 3.2, there is always a sequence $\{z_k\}$
in ${\s X}^{**}$ which converges weak-star to $0$ and
which satisfies
$
||z_k|| = 1 = \text{\rm dist} (z_k,
\text {\rm sp} \{z_1, \ldots , z_{k-1} \} ).
$
Thus we obtain the following.

\corollary 5.2
{\it If $\s X$ is an infinite-dimensional Banach space, then there
are analytic disks in the fiber ${\s M}_0$ over $0$ of the
spectrum $\s M$ of $\Hi (B)$.
In fact, there is an analytic embedding of
$(\beta (N) \backslash N) \times \Delta $ into ${\s M}_0$.}

\vskip .3in

{\bf 6. Embedding infinite dimensional analytic structure
in $\s M_0.\ $}
The question now arises as to what sort of analytic objects are to be
found in the fibers of the spectrum of $\Hi (B)$.  When is it possible
to inject analytically the unit ball $B$ of $\s X$ into $\s M _0$?  What
are the natural analytic objects appearing in $\s M_0$?  Can one learn
something about $\s X$ by peering into $\s M_0$?

If the unit ball of some infinite dimensional Banach space injects
analytically in $\s M_0$, then so does the (countable dimensional)
infinite polydisk $\Delta ^\infty$, the unit ball of $\ell ^\infty$.
This is because $\ell ^\infty$ can be mapped injectively into any
infinite dimensional Banach space.  In the other direction, any
separable Banach space maps injectively into $\ell ^\infty$, so if one
can inject $\Delta ^\infty$ analytically into $\s M_0$, then one can
inject the unit ball of any separable Banach space analytically into
$\s M_0$.

Conceptually the simplest way to inject $\Delta ^\infty$ is by shifting
along a basic sequence and passing to a limit. In Theorem 6.1 we give
one result with this method of proof.
The reader should note that all the results of this section through
Theorem 6.5 apply to the spaces
$\ell ^p$ and $L^p[0,1]$ for $1 < p < \infty$.

\theorem 6.1
{\it Suppose $\s X$ has a normalized basis $\{x_j\}$ which is shrinking,
i.e., whose associated coefficient
functionals $\{L_j\}$ have linear span dense in $\s X^*$.
Suppose furthermore that there is an integer $N \geq 1$ such that
$$
\sum |L_j(x)|^N < \infty
$$
for all $x = \sum L_j(x) x_j$ in $\s X.$
Then there is an analytic injection of the (countable dimensional)
infinite polydisk $\Delta ^\infty$ into the fiber $\s M_0.$ }

\proof
Define a map $T_k$ of $\ell ^ \infty$ into $\s X$ by
$$
T_k(y) = \sum _{n=1}^{\infty} 2^{-n}y_nx_{n+k}, \qquad y \in
\ell ^ \infty .
$$
Then $T_k$ injects $\Delta ^\infty$ into $B$.  Since $\s M$ is compact,
the maps from $B$ to $\s M$ form a compact set in the topology of
pointwise convergence.  Let $T$ be any map which is adherent to the
sequence $T_k$ as $k \rightarrow \infty$.
Evidently $f \circ T$ is analytic on $B$ for all $f \in
\Hi (B)$, as it is a cluster point of the sequence $f \circ T_k$ of
bounded analytic functions on $B$.  Thus $T$ is an analytic mapping of
$\Delta ^\infty$ into $\s M$.

If $L \in \s X^*$ is a finite linear combination of the coordinate
functionals $L_j$, then $L \circ T_k$ converges pointwise to $0$ on $B$,
so that $L \circ T = 0$.  Since these finite linear combinations are
dense in $\s X ^*$, every $L \in \s X ^*$ vanishes on the image of $T$,
and the image of $T$ is contained in the fiber $\s M_0$.

Suppose $u$ and $v$ are distinct points of $\Delta ^ \infty$ with $Tu
= Tv$.  By the principle of uniform boundedness, there is $C > 0$
such that $\sum |L_j(x)|^N \leq C$ for all $x \in B$.  Choose a net
$k_{\alpha} \rightarrow \infty$ such that $T_{k_\alpha}$ converges
pointwise on $\Delta ^ \infty$ to $T$.  For
$\lambda \in {\text {\bf C}}$ with $|\lambda| = 1$, define an
$N$-homogeneous analytic function $f _\lambda$ on $\s X$ by
$$
f _\lambda (x) = \sum _{j=1}^\infty \lambda ^j L_j (x) ^N \ ,
\qquad x \in \s X.
$$
Evidently $|f _\lambda| \leq C$ on $B$, and
$$
f_\lambda (T_k(y))
= \lambda ^ k \sum _{j=1} ^ \infty \lambda ^j ( 2 ^{-j} y _j ) ^ N
= \lambda ^k f _ \lambda ( T _0 ( y ) ),
\qquad y \in \Delta ^ \infty , k \geq 0.
$$
We may assume $\chi (\lambda ) = \lim  \lambda ^{k _ \alpha}$
exists for all $|\lambda | = 1$. From the formula above
we obtain
$ f _ \lambda ( T ( y ) ) = \chi (\lambda) f _ \lambda ( T _0 ( y ) )$
for all $y \in \Delta ^\infty$. Since $|\chi (\lambda )| = 1$ and
$Tu = Tv$, we have $f_ \lambda (T_0(u))=f_ \lambda (T_0(v))$
for all $\lambda$.
Equating coefficients in the power series expansion, we obtain
$u _j ^N = v _j ^N$ for all $j \ge 1$.
Replacing $N$ by $N+1$ in this argument, we also obtain
$u_j^{N+1} = v_j^{N+1}$ for all $n \ge 1$. It follows that $u = v$,
and $T$ is one-to-one.
\endproof

For subsequent embedding theorems,
we will use the {\it Gleason metric} $\rho$ on $\s M$, defined so
that $\rho (\varphi,\psi)$ is the supremum of $|f(\varphi) - f(\psi)|$
over all $f \in \Hi (B)$ satisfying $||f|| \leq 1$.  (See Chapter VI
of [Ga].) The Gleason metric
on $B$, regarded as a subset of $\s M$, is equivalent to the norm
metric.  In fact, it is a simple consequence of the Schwarz lemma
that for any fixed $r < 1$, the Gleason metric is uniformly
equivalent to the norm metric on the ball $B_r = \{ ||x|| \leq r\}$
via a bi-Lipschitz map.

A more natural problem, and a more difficult problem, than the problem
of simply embedding analytic structure is to embed
analytic structure in $\s M _0$ via an analytic injection which is
uniformly bicontinuous with respect to the appropriate Gleason metrics.
To obtain such an injection, we will look for limits along a sequence of
subspaces.

Fix $p > 1$.  A sequence $\{ {\s W}_n \}$ of
subspaces of $\s X^*$ satisfies an {\it upper $p$-estimate} if
for some constant $C$,
$$
||\sum L_n|| \leq C\{ \sum ||L_n||^p \}^{1/p}
$$
whenever $L_n \in {\s W}_n  , n \geq 1$. Let $q = p/(p-1)$ be the conjugate
index to $p$. Then for any integer $N \geq q$ and any functionals
$L_n \in {\s W}_n$ satisfying $||L_n|| \leq 1$, the series
$$
f = \sum _{n=1} ^{\infty} L_n^N
$$
converges to an analytic $N$-homogeneous function on $\s X$ which
satisfies
$$
|f(x)| \leq C^q, \qquad ||x|| \leq 1.
$$
Indeed, $|f(x)|$ is bounded by
$$
\sum |L_n(x)|^N \leq \sum |L_n(x)|^q
\leq \sup_{\sum |a_n|^p \le 1} |\sum a_n L_n(x)|^q
\leq C^q.
$$

Let $\s U$ be any free ultrafilter on the positive integers.
Recall that the {\it $\s U$-ultraproduct} of a sequence of Banach
spaces is the quotient space obtained from their
$\ell^\infty$-direct sum by dividing out by the subspace of
sequences whose norms have limit $0$ along the ultrafilter.
The norm of an element $\tilde x$ represented by a sequence
$(x_1,x_2,...)$ is given by
$||{\tilde x}|| = \lim_{\ss{{\s U}}} ||x_n||$.
As an example, the limit of the scalar products in $\ell ^2_n$
determine a scalar product which makes $\lim_{\ss{{\s U}}}
\ell ^2_n$ into a Hilbert space of uncountable dimension.
For background on ultraproducts see [He].

Recall that a subspace $\s W$ of ${\s X} ^*$ is said to
$c-${\it norm} the subspace $\s E$ of $\s X$ if
the supremum of $|L(x)|$ over $L \in \s W$,
$||L|| \leq c,$ is at least $||x||$ for every $x \in \s E$.
One defines similarly what it means for $\s E$ to $c-${\it norm}
$\s W$.

\lemma 6.2
{\it Fix $p > 1$ and $C,c \geq 1$.
Let $\{{\s E}_n\} $ be a sequence of subspaces of ${\s X}$ and
$\{{\s W}_n\}$ a sequence of subspaces of ${\s X^*}$ such that
${\s W}_n$ c-norms ${\s E}_n$ for each $n$,
while ${\s W}_n$ is orthogonal to ${\s E}_k$ for all $k \neq n$.
Assume that the sequence $\{ {\s W}_n \}$ satisfies an upper
$p$-estimate, while the sequence $\{ {\s E}_n \}$
tends weakly to $0$ (i.e., if $x_n\in {\s E}_n$ is uniformly
bounded, then $x_n\rightarrow 0$ weakly).
Then for any free ultrafilter $\s U$ on the positive integers,
there is a natural analytic injection of
the unit ball $\tilde B$ of the ${\s U}$-ultraproduct
$ \tilde {\s E}$ of $\{{\s E}_n\} $ into
the fiber ${\s M_0}$ of the spectrum of $\Hi (B)$
which is uniformly bicontinuous, from the norm metric
of any ball $r \tilde B , r < 1,$ to the Gleason metric of $\s M$. }

\proof
Suppose $\tilde x \in \tilde {\s E}$ satisfies $||\tilde x|| < 1$.
Then $\tilde x$ is represented by
$(x_1,x_2,\dots)$ with $x_n \in {\s E}_n$ satisfying $||x_n|| < 1$.
Since $\s M$ is compact,
$\lim_{\ss{{\s U}}} x_n $ exists
in ${\s M}$. If $(y_1,y_2,\dots)$ also represents $\tilde x$,
then $\lim_{\ss{{\s U}}} ||x_n - y_n|| = 0$, so that
$\lim_{\ss{{\s U}}} \rho (x_n,y_n) = 0$, and consequently
$\lim_{\ss{{\s U}}} x_n = \lim_{\ss{{\s U}}} y_n.$
We denote this common limit by $S(\tilde x)$.
As before one sees that $S$ is an analytic
map, from $\tilde B$ into ${\s M}$. Any such analytic map is nonexpansive,
from the Gleason metric of $\tilde B$ determined by $\Hi (\tilde B)$ to
the Gleason metric of $\s M$.
Since $\{{\s E}_n\}$
tends weakly to $0$, $S$ maps into ${\s M_0}$.

It remains to check that $S$ is one-to-one with uniformly continuous
inverse.  Take any $\tilde x , \tilde y$ in $\tilde B$ with
$||\tilde x - \tilde y || > \delta > 0$ and choose $U \in {\s U}$ so
that $||x_n - y_n || > \delta $ for all $n\in U$.  For $n \in U$ pick
$L_n$ in the unit ball of ${\s W}_n$ so that $|L_n(x_n)-L_n(y_n)| >
\delta /c$.  Either for $N=[q]+1$ or $N=[q]+2$ the set $V =
\{n\in U \colon \ |\ [L_n(x_n)]^N - [L_n(y_n)]^N | > \tau \}$ (where
$\tau > 0$ depends on $c$ and $p$ but not on $n$) is in the ultrafilter
${\s U}$.  Set $f=\sum_{\ss{n \in V}} L_n^N$.  Then by the earlier
remark, $f$ is in $\Hi (B)$ and $||f|| \leq C^q$.  But for each $n \in
V$, $|\ [f(x_n)] - [f(y_n)]\ | = |\ [L_n(x_n)]^N - [L_n(y_n)]^N | >
\tau$, so $|f(\tilde x)-f(\tilde y)| \geq \tau$ and
$C^q \rho(S(\tilde x),S(\tilde y)) \geq \tau$.
\endproof

A Banach space $\s Y$ is {\it finitely representable} in $\s X$ if
for any $c > 1$, any finite dimensional subspace of $\s Y$ is
isomorphic to a subspace of $\s X$ via an isomorphism $T$ which
satisfies $||T|| < c$ and $||T^{-1}|| < c$. The Banach space $\s X$
is {\it superreflexive} if any Banach space finitely representable
in $\s X$ is reflexive. According to a theorem of P. Enflo, a Banach
space is superreflexive if and only if it is uniformly convexifiable,
and this occurs if and only if its dual space is so.
For a discussion of this circle of ideas, see pp. 86-87 of [Di1].

\theorem 6.3
{\it If ${\s X}$ is a superreflexive Banach space, then the unit ball
of a non-separable Hilbert space injects into the fiber
${\s M_0}$ via an analytic map which is uniformly bicontinuous from the
norm metric of the unit ball of the Hilbert space to the Gleason metric
of its image in $\s M$.}

\proof
With Dvoretzky's theorem and a standard gliding hump argument,
it is possible to construct for any infinite dimensional Banach
space $\s X$ sequences of finite dimensional
subspaces ${\s E}_n$ of $\s X$ and ${\s W}_n$
of $\s X^*$ such that the dimension of ${\s E}_n$ tends to $\infty$,
each ${\s W}_n$ c-norms ${\s E}_n$,
${\s W}_n$ is orthogonal to ${\s E}_k$ for all $k \neq n$, and
the Banach-Mazur distance $d({\s E}_n, \ell^2_n)$ tends to $1$.
The inductive step for constructing $\{{\s E}_n\}$ and $\{{\s W}_n\}$
goes as follows.  Suppose that the first $n$ subspaces in each
sequence have been constructed.
Choose a finite dimensional subspace
${\s E} \supset \cup_{k=1}^n {\s E}_k$ of ${\s X}$ which almost
precisely norms the span of $\cup_{k=1}^n {\s W}_k$ and
choose a finite dimensional subspace
${\s W} \supset \cup_{k=1}^n {\s W}_k$ of ${\s X^*}$ which almost
precisely norms  ${\s E}$. A simple argument based on the Hahn-Banach
theorem shows that
${\s E}^{\perp}$ almost 2-norms ${\s W}_{\perp}$. By Dvoretzky's
theorem we may choose a subspace
${\s E}_{n+1}$ of ${\s W}_{\perp}$ of arbitrarily large dimension
which is almost isometric to $\ell^2_{n+1}$. Then choose
${\s W}_{n+1} \subset {\s E}^{\perp}$ finite dimensional to
almost 2-norm ${\s E}_{n+1}$.

If $\s W$ norms $\s E$ above with enough precision, a standard argument
([Di2], pp 38-39) shows that any sequence $\{ x_n \}$ with $x_n$ a
nonzero vector in ${\s E}_n$ forms a basic sequence.  Since any
normalized basic sequence in a reflexive Banach space tends weakly to
$0$ (because the coefficient functionals vanish on any weak adherent
point), the sequence of subspaces ${\s E}_n$ tends weakly to $0$.

If $\s E$ norms $\s W$ above with enough precision, the same standard
argument shows that any sequence $\{ L_n \}$ with $L_n$ a unit vector in
${\s W}_n$ forms a basic sequence, and furthermore the basis constants
are uniformly bounded, independent of the normalized basic sequence.
Thus the Gurarii-Gurarii-James theorem shows that
any such sequence satisfies an upper $p$-estimate for some $p > 1$, and
moreover by Theorem 2 of [Ja] (or by the discussion in [Di2]),
the index $p$ can be chosen independent of the
normalized basic sequence $\{ L_n \}$.

The hypotheses of Lemma 6.2 are now satisfied.
Since $d({\s E}_n, \ell^2_n) \rightarrow 1$, any nontrivial
ultraproduct of $\{{\s E}_n\}$ is a Hilbert space of uncountable
dimension. If we restrict the map of Lemma 6.2 to a smaller
ball, we obtain a uniformly bicontinuous injection.
\endproof

The same method of proof establishes the following.

\theorem 6.4
{\it If ${\s X}$ is superreflexive and
$\ell^p$ is finitely representable in ${\s X}$ for
some $1\leq p < \infty$, then the unit balls of
$L^p[0,1]$ and of $\ell^p_{\aleph}$
inject into ${\s M_0}$ via analytic maps which are uniformly
bicontinuous from the norm of the unit ball to the Gleason
metric of $\s M$.}

\proof
If $\ell^p$ is finitely representable in ${\s X}$ then it
is also finitely representable in every finite
codimensional subspace of ${\s X}$.  The proof of
Theorem 6.3 then shows that every ultraproduct of
$\{\ell^p_n\}$ injects into ${\s M_0}$ via
an analytic uniformly bicontinuous map, while $L^p[0,1]$ and
$\ell^p_{\aleph}$ inject isometrically into every nontrivial
ultraproduct of $\{\ell^p_n\}$ (see [Hn], or [He]).
(To see that $L^p[0,1]$ injects, consider maps to $\ell ^p_n$
obtained by averaging functions over intervals of length $1/n$.
To see that $\ell^p_{\aleph}$ injects, consider for each
$0 < t \leq 1$ the $\s U$-limit $u_t$ of the sequence $v_n \in
\ell^p_n$ with only one nonzero entry, a $1$ in the coordinate
position $[tn]$, the integral part of $tn$.)
\endproof

In the case that $\s X$ is $L^p[0,1]$ or $\ell ^p$ for
$1 < p < \infty$, we can embed all of $\s M$ into $\s M_0$.
In fact, we have the following.

\theorem 6.5
{\it Suppose $1 < p < \infty$, and suppose $\s X$ is isometric to
an infinite $\ell ^p$-direct sum of itself. Then there is a
homeomorphism $\Phi$ of $\s M$ onto a compact subset of $\s M _0$
satisfying $f \circ \Phi \in \Hi (B)$ for all $f \in \Hi (B)$,
which is uniformly bicontinuous with respect to the Gleason metric.}

\proof
For this proof, denote by ${\s X} _j$ a copy of $\s X$, let $\tilde
{\s X}$ denote the $\ell ^p$-direct sum of the ${\s X} _j$'s, and
let $\tilde B$ denote the unit ball of $\tilde {\s X}$, which is
by hypothesis isometric to $B$. Fix a free ultrafilter $\s U$ on
the integers, and define an operator $\Lambda$ from $\Hi (\tilde B)$
to $\Hi (B)$ so that $(\Lambda f)(x)$
is the pointwise limit along $\s U$ of the sequence of
functions $f(0,\ldots,0,x,0,\ldots )$, where the $x$ appears in the
$n$th coordinate. Evidently $\Lambda$ is linear and continuous and
multiplicative, and $||\Lambda ||=1=\Lambda (1)$. The dual map
$\Phi$ is a continuous map from ${\s M} (B)$ into ${\s M}(\tilde B)$
satisfying $f \circ \Phi \in \Hi (B)$ for all $f \in \Hi (\tilde B)$.
Since ${\tilde {\s X}}^*$ is the $\ell ^q$-direct sum of the duals of
the ${\s X} _j$'s, it is clear that $\Lambda (f) = 0$ for all
$f \in {\tilde {\s X}}^*$. Hence the image of ${\s M}(B)$ under $\Phi$
is contained in ${\s M} _0(\tilde B)$.

Let $p_0$ be the least integer satisfying $p_0 \geq p$.  Suppose $g \in
\Hi (B)$ has order at least $p_0$ at the origin, that is, the terms of
the Taylor series of $g$ at $0$ of order less than $p$ are all zero.
Define an analytic function $G$ on $\tilde {\s X}$ by
$$
G(x_1\,,x_2\,,\ldots) = \sum_{j=1}^{\infty} g(x_j), \qquad
(x_1\,,x_2\,,\ldots) \in \tilde {\s X}.
$$
The estimate $|g(x_j)| \leq ||g||_B ||x_j||^p$ leads to $||G||_{\tilde
B} \leq ||g||_B$.  (In fact, equality holds here.)  Evidently $\Lambda
(G) = g$.  Since such $g$'s separate the points of ${\s M}(B)$, the map
$\Phi$ is one-to-one, hence a homeomorphism onto its image.

The map $\Phi$ is non-expanding with respect to the Gleason
metric. Let $\varphi, \psi \in {\s M}(B)$. Suppose $h \in \Hi (B)$
satisfies $||h|| \leq 1$ and $|h(\varphi) - h(\psi )| = d$.
Then $2^{-p_0}(h - h(0))^{p_0}$ and $2^{-p_0-1}(h - h(0))^{p_0+1}$
vanish to order at least $p$ at $0$ and have norms at most $1$,
and either one or the other of these functions satisfies
$|g(\varphi) - g(\psi )| \geq c d^{p+2}$, where $c$ depends only on $p$.
Hence the corresponding $G$ satisfies
$|G(\Phi (\varphi)) - G(\Phi(\psi ))| \geq c d^{p+2}$. This yields the
uniform bicontinuity of $\Phi$ with respect to the Gleason metric.
\endproof

Theorem 6.1 does not apply to the Banach space $c_0$.  By making use
of the automorphisms of the unit ball of $c_0$, we can still inject
infinite dimensional analytic objects in $\s M_0$.

\theorem 6.6
{\it Let $\s X$ be the Banach space $c_0$ of null sequences.  There is
an analytic injection of the unit ball $B^{**}$ of $\ell ^ \infty $ into
the fiber ${\s M_0}$ which is an isometry from the Gleason metric of
$B^{**}$ to the Gleason metric of $\s M$. }

\proof
Choose  $0 < \alpha_k < 1$  such that  $\alpha_k$
increases to 1 and  $\Sigma(1 - \alpha_k) < \infty$, and choose integers
$m_k$  such that  $m_{k+1} > m_k + k$.  Define
$\Phi_k : B^{**} \rightarrow B$ by
$$
\Phi_k (w) = \bigg(0,\ldots,0, \frac{\alpha_k - w_1}{1 -
\alpha_k w_1},\ldots,\frac{\alpha_k - w_k}{1 - \alpha_k w_k} , 0,
\ldots\bigg), \qquad w \in B^{**},
$$
where the first block consists of $m_k$ zeros.
As in the proof of Lemma 6.2,
the sequence  $\{\Phi_k\}$  has a pointwise limit along
any fixed free ultrafilter $\s U$,
which determines an analytic map $S$ of $B^{**}$ into $\s M$.
The map $S$ is nonexpansive, from the Gleason metric of $B^{**}$
to the Gleason metric of $\s M$.
Evidently the image of $S$ annihilates $\ell ^1$, and so is
contained in the fiber $\s M _0$ over $0$.

Next observe that the Gleason metric $\rho$ of $B^{**}$ is expressed
in terms of the Gleason metric $\rho_{\ss{\Delta}}$ of the open unit disk
$\Delta$ by
$$
\rho (z,w) = \sup _j \rho_{\ss{\Delta}}(z_j\, , w_j), \qquad z,w \in B^{**}.
\tag {6.1}
$$
The estimate $\rho_ {\ss{\Delta}} (z_j\, ,w_j) \leq \rho (z,w)$ is
obtained by considering functions which depend only on the $j$th
coordinate of $B^{**}$. Taking the supremum over $j$, we obtain
the inequality $\geq$ in (6.1).
To establish the reverse inequality, we can
by composing each coordinate variable with an
appropriate linear fractional transformation assume that $z = 0$.
Consideration of the
analytic disk
$$
\lambda \rightarrow \frac {1}{||w||} (\lambda w_1\, ,\lambda w_2 \, ,
\cdots ), \qquad |\lambda |<1,
$$
yields then easily the estimate $\rho (0,w) \leq \rho_{\ss{\Delta}}
(0,||w||)$,
which coincides with the supremum of $\rho _{\ss{\Delta}} (0,w_j)$. This
establishes (6.1).

For $N \geq i \geq 1$ define
$$
g_{i,N}(z) = \prod_{j=N}^\infty \dfrac{\alpha_j - z_{i+m_j}}{1
 - \alpha_j z_{i+m_j}},
\qquad z \in B^{**}.
$$
Our conditions on the $\alpha _j$'s guarantee that the product converges,
and  $g_{i,N} \in H^\infty(B)$.
We compute that
$$
g_{i,N} (\Phi _k(z)) = z_i
\prod_{j \geq N,j\neq k} \ \alpha_j, \qquad k \geq N.
$$
Taking the limit as $k$ tends to $\infty$ through $\s U$, we obtain
$$
g_{i,N} (S(z)) = z_i
\prod_{j \geq N} \ \alpha_j.
$$
For $z,w \in B^{**}$, choose $h \in \Hi (\Delta)$ such that
$|h| \leq 1$ and $\rho _{\ss{\Delta}} (z_i\, ,w_i) = |h(z_i) - h(w_i)|.$
Then $h \circ g_{i,N}$ belongs to $\Hi (B^{**})$, and
$(h \circ g_{i,N})(S(z)) - (h \circ g_{i,N})(S(w))$ tends to
$\rho _{\ss{\Delta}} (z_i\, ,w_i)$ as $N \rightarrow \infty$.
It follows that
$$
\rho(S(z),S(w)) \geq \rho _{\ss{\Delta}} (z_i\, ,w_i), \qquad i \geq 1.
$$
Taking the supremum over $i$ and using (6.1), we obtain
$$
\rho (S(z),S(w)) \geq \rho (z,w), \qquad z,w \in B^{**}.
$$
Since $S$ is nonexpansive, equality holds here, and the proof is
complete.
\endproof

The analytic structure produced in Theorem 6.6 is by no means
maximal. In fact, it is possible to fit $S$ into a
family $\zeta \rightarrow S_\zeta$ of analytic injections,
depending analytically
on the parameter $\zeta \in B^{**}$, so that the image of
$B^{**}$
under $S_\zeta$ is contained in the fiber over $\zeta.$
To do this one simply introduces an analytic parameter $\zeta$
in the definition of
$\Phi _k$, setting
$$
\Phi_k (\zeta ,w) = \bigg(\zeta _1\ ,\ldots,\zeta _k\ ,
0,\ldots,0, \frac{\alpha_k - w_1}{1 -
\alpha_k w_1},\ldots,\frac{\alpha_k - w_k}{1 - \alpha_k w_k} , 0,
\ldots\bigg),\qquad \zeta ,w \in B^{**},
$$
and one passes to the limit along the free ultrafilter as above,
arriving thereby at the following result.

\theorem 6.7  {\it Let  ${\s X} = c_0$, and let  $B$  and
$B^{**}$  be the open unit balls of  $c_0$  and  $\ell^\infty$
respectively.  Then there is a map
$$
\Phi : B^{**} \times B^{**} \rightarrow {\s M}
$$
such that:
\roster
\item"{\rm (i) \ }" $\Phi$ is continuous,
with respect to the norm topology on $B^{**}$.
\item"{\rm (ii) \ }" $f\,{}_{{}^\circ} \Phi$ is
analytic on $B^{**} \times B^{**}$ for all
$f \in H^\infty(B)$.
\item"{\rm (iii) \ }" $\Phi(z,w) \in {\s M}_z$ for all
$z, w \in B^{**}$.
\item"{\rm (iv) \ }" The image of $\Phi$ is disjoint
from $B^{**}$.
\item"{\rm (v) \ }" $\Phi$ is one-to-one.
\endroster
}

\vskip .4in

{\bf 7. Other algebras related to $\Hi(B)$.} \ To preserve consistency
of point of view, we have limited our investigation to the algebra
$\Hi(B)$.
Nonetheless, there are other natural algebras of analytic functions
on $B$,
and analytic structure in ${\s M}$
is related to analytic structure in the maximal ideal spaces
of these other algebras.

First of all, there is the algebra $H _b ({\s X})$ (see [ACG])
which consists of all entire functions on $\s X$
that are bounded on bounded sets and which has
spectrum ${\s M }_b({\s X})$.
The restrictions of these entire functions to $B$
generate a uniform algebra denoted $\Hi _{uc} (B)$,
the maximal ideal space of which can be identified with
a compact subset of ${\s M }_b$.
Lastly, there is $A(B)$, the infinite dimensional analogue of
the disk algebra, which is the uniformly closed algebra generated
on $B$ by ${\s X} ^*$.
It is easily seen that its maximal ideal space
is  just $\bar B^{**}$.

Dual to the inclusions
$A(B) \subset \Hi _{uc} (B) \subset \Hi (B)$,
there exist induced maps between the corresponding
maximal ideal spaces:
$$
{\s M} = {\s M} _{\Hi (B)}
\overset \pi _1 \to \longrightarrow 
{\s M} _{\Hi _{uc} (B)}
\overset \pi _2 \to \longrightarrow 
{\s M} _{A(B)} = \bar B^{**}.
$$
Clearly, $\pi _2 \circ \pi _1 : {\s M} \rightarrow \bar B ^{**}$
is just the projection used elsewhere in this paper.
It is important to look for analytic structure in the fibers of both
$\pi _1$ and $\pi _2$.
A careful examination of the results given earlier in this paper
shows that sometimes analytic structure lies in fibers of
$\pi _1$, sometimes in those of $\pi _2$.

For example,
Corollary 5.2 admits the following refinement.

\theorem 7.1
{\it
Let ${\s X}$ be an infinite dimensional Banach space.
Then,
$\pi _1 ^{-1} (\zeta)$ contains an analytic disk for some
$\zeta \in \pi _2 ^{-1} (0)$.
}

\proof
As in the proof to Theorem 5.1, there exist a disk
$D \subset {\text {\bf C}}$ and analytic maps
$z _k : D \rightarrow B^{**}$
so that, for an appropriate subnet $\{ k _\alpha \}$,
$z = \lim _{\alpha} z _{k _\alpha}$ defines a one-to-one,
analytic map
$z : D \rightarrow {\s M}_0$.
The one piece of new information we provide here is that
$f \circ z$ is constant on $D$ for each $f \in \Hi _{uc} (B)$.
For that purpose, it suffices to assume $f$ is homogeneous since
the linear span of such functions is dense in $\Hi _{uc} (B)$.
In [DG] it is shown that the canonical extension of a homogeneous
$f$ is homogeneous and hence uniformly continuous with respect
to the norm topology on $B^{**}$.
By construction, the maps $z _k$ also satisfy
$|| z _k (\lambda) - z _k(0) || < c _k$
for $\lambda \in D$ where $c _k  \rightarrow 0$;
so uniform continuity implies $f \circ z$ is constant.
Hence, we conclude that $ \pi _1 \circ z$ maps $D$ to a single point
$\zeta \in {\s M} _{\Hi _{uc} (B)}$ where $ \pi _2 ( \zeta ) = 0$.
\endproof

In a very different direction, the results of Section 6 through
Theorem 6.5 establish the existence of analytic maps into
${\s M}_0$ by studying entire functions on ${\s X}$.
For each of these theorems, we could simply replace
${\s M}_0$ with $\pi _2 ^{-1} (0)$; the proofs would remain unchanged.
Consequently, these results give information about analytic structure
in fibers of $\pi _2$ and shed no additional light on analyticity
in the fibers of $\pi _1$.
Moreover, use of entire functions in section 6 means that most
of the arguments apply to the algebra
$H _b ({\s X})$ of entire functions
and its spectrum ${\s M} _b({\s X})$.
For example, the proof to Theorem 6.3 shows the following

\theorem 7.2
{\it If ${\s X}$ is a superreflexive Banach space, then
a non-separable Hilbert space injects into
${\s M} _b({\s X})$
via an analytic map.}

Finally, for ${\s X} = c _0$, it is known (see [ACG])
that $\Hi _{uc} (B) = A(B)$.
So, the fibers of $\pi _2$ are trivial,
and the analytic structure described in Theorem 6.6
actually lies in $\pi _1 ^{-1} (0)$.

In conclusion we mention that 
we know of no nontrivial restriction on the
analytic structure in ${\s M}_0$. There are a number of
test questions which can be raised.
For instance,
does there always exists infinite dimensional analytic
structure in ${\s M}_0$, or better yet, in
$\pi _1 ^{-1} (0)$?
More specifically, does the unit ball $\Delta ^ \infty$
of $\ell ^ \infty$ inject analytically into ${\s M}_0$ via
a uniformly bicontinuous mapping when $\s X = \ell ^2$? A
negative answer to the latter would suggest a real connection
between the function theory and the geometry of the Banach space.

\newpage

\centerline{REFERENCES}

\roster

\item"{[ACG] \ }"
R.~M. Aron, B.~J. Cole and T.~W. Gamelin, Spectra of algebras of
analytic functions on a Banach space, J. reine angew. Math. 415
(1991).

\item"{[BFT] \ }"
J. Bourgain, D.~H. Fremlin and M. Talagrand, Pointwise compact sets
of Baire-measurable functions, Amer. J. Math. 100 (1978), 845-886.

\item"{[DG] \ }"
A.~M. Davie and  T.~W. Gamelin, A theorem on polynomial-star
approximation, Proc. Amer. Math. Soc. 106 (1989), 351-356.

\item"{[Di1] \ }"
J. Diestel, {\it Geometry of Banach Spaces - Selected Topics},
Lecture Notes in Math., Vol. 485, Springer-Verlag, 1975.

\item"{[Di2] \ }"
J. Diestel, {\it Sequences and Series in Banach Spaces}, Springer-Verlag, 1984.

\item"{[Ga] \ }"
T.~W. Gamelin, {\it Uniform Algebras}, Second Edition, Chelsea Press,
1984.

\item"{[Gl] \ }"
J. Globevnik, Boundaries for polydisc algebras in infinite
dimensions, Math. Proc. Camb. Phil. Soc. 85 (1979), 291-303.

\item"{[He] \ }"
S. Heinrich, Ultraproducts in Banach space theory,
J. reine angew. Math. 313 (1980), 72-104.

\item"{[Hn] \ }"
C. Ward Henson, Non standard hulls of Banach spaces, Israel J.
Math. 25 (1976), 109-144.

\item"{[HJ] \ }"
J. Hagler and W.~B. Johnson, On Banach spaces whose dual balls
are not weak$^*$ sequentially compact, Israel J. Math. 28 (1977),
325-330.

\item"{[Ho] \ }"
K. Hoffman, Bounded analytic functions and Gleason parts, Annals of
Math. 86 (1967), 74-111.

\item"{[Ja] \ }"
R.C. James, Super-reflexive spaces with bases, Pac. J. Math.
41 (1972), 409-419.

\item"{[Jo] \ }"
B. Josefson, Weak sequential convergence in the dual of a Banach space
does not imply norm convergence, Arkiv f\"or mat. 13 (1975), 79-89.

\item"{[KKM] \ }"
M.A. Krasnoselskii, M.G. Krein and D.P. Milman, On the deficiency
indices of linear operators in Banach spaces and some geometrical
questions, Sbornik Trudov Inst. Akad. Nauk Ukr. SSR
11 (1948), 97-112.

\item"{[LT] \ }"
J. Lindenstrauss and L. Tzafriri, {\it Classical
Banach Spaces I, Sequence Spaces},
Springer-Verlag, 1977.

\item"{[MS] \ }"
V.~D. Milman and G. Schechtman, {\it Asymptotic Theory of
Finite Dimensional Normed Spaces,} Lecture Notes in Math., Vol. 1200,
Springer-Verlag, 1986.

\item"{[Mu] \ }"
J. Mujica, {\it Complex Analysis in Banach Spaces,} Math. Studies 120,
North Holland, 1986.

\item"{[Ni] \ }"
A. Nissenzweig, $w^*$ sequential convergence, Israel J. Math. 22
(1975), 266-272.

\item"{[ORS] \ }"
E.~W. Odell, H.~P. Rosenthal, and T. Schlumprecht,
On weakly null FDD's in Banach spaces, to appear.

\item"{[Pe] \ }"
A. Pe\l czy\'nski, On Banach spaces containing $L^1(\mu )$, Studia Math.
30 (1968), 231-246.

\endroster

\end